\documentclass[a4paper,12pt,twoside,leqno]{article} 
\usepackage[english]{babel} 
\usepackage{amsmath,amssymb,amsthm,amsfonts} 
\usepackage[dvips]{graphicx}
\usepackage[all]{xy}

\usepackage{hyperref}
\hypersetup{
  pdftitle={Steinitz classes of tamely ramified Galois extensions of algebraic number fields},
  pdfauthor={Alessandro Cobbe},
  pdfsubject={Steinitz classes},
  pdfkeywords={Steinitz classes}}

\usepackage{fancyhdr}
\newcommand{\Q}{\mathbb Q}
\newcommand{\N}{\mathbb N}
\newcommand{\Z}{\mathbb Z}

\renewcommand{\P}{\mathfrak P}
\newcommand{\qq}{\mathfrak Q}
\newcommand{\p}{\mathfrak p}
\newcommand{\q}{\mathfrak q}
\newcommand{\gal}{\mathrm{Gal}}
\renewcommand{\epsilon}{\varepsilon}
\newcommand{\disc}{\mathrm{d}}
\newcommand{\norm}{\mathrm{N}}

\newcommand{\st}{\mathrm{st}}

\newcommand{\rt}{\mathrm{R}_t}
\newcommand{\cl}{\mathrm{Cl}}
\newcommand{\aut}{\mathrm{Aut}}
\newcommand{\ind}{\mathrm{ind}}
\newcommand{\res}{\mathrm{res}}
\renewcommand{\hom}{\mathrm{Hom}}
\newcommand{\G}{\mathcal G}
\newcommand{\oo}{\mathcal O}

\newcommand{\comment}[1]{}

\newtheorem{teo}{Theorem}[section]
\newtheorem{lemma}[teo]{Lemma}

\newtheorem{prop}[teo]{Proposition}
\newtheorem{defn}[teo]{Definition}

\title{Steinitz classes of tamely ramified Galois extensions of algebraic number fields}
\author{Alessandro Cobbe}

\begin{document}
\maketitle

\markboth{Abstract}{Abstract}
\begin{section}*{Abstract}
\addcontentsline{toc}{section}{Abstract}
The Steinitz class of a number field extension $K/k$ is an ideal class in the ring of integers $\oo_k$ of $k$, which, together with the degree $[K:k]$ of the extension determines the $\oo_k$-module structure of $\oo_K$. We call $\rt(k,G)$ the classes which are Steinitz classes of a tamely ramified $G$-extension of $k$. We will say that those classes are realizable for the group $G$; it is conjectured that the set of realizable classes is always a group.

We define $A'$-groups inductively, starting by abelian groups and then considering semidirect products of $A'$-groups with abelian groups of relatively prime order and direct products of two $A'$-groups. Our main result is that the conjecture about realizable Steinitz classes for tame extensions is true for $A'$-groups of odd order; this covers many cases not previously known. Further we use the same techniques to determine $\rt(k,D_n)$ for any odd integer $n$.

In contrast with many other papers on the subject, we systematically use class field theory (instead of Kummer theory and cyclotomic descent).
\end{section}

\markboth{Introduction}{Introduction}
\begin{section}*{Introduction}
\addcontentsline{toc}{section}{Introduction}

Let $K/k$ be an extension of number fields and let $\oo_K$ and $\oo_k$ be their rings of integers. By Theorem 1.13 in \cite{Narkiewicz} we know that
\[\oo_K\cong \oo_k^{[K:k]-1}\oplus I\]
where $I$ is an ideal of $\oo_k$. By Theorem 1.14 in \cite{Narkiewicz} the $\oo_k$-module structure of $\oo_K$ is determined by $[K:k]$ and the ideal class of $I$. This class is called the \emph{Steinitz class} of $K/k$ and we will indicate it by $\st(K/k)$. Let $k$ be a number field and $G$ a finite group, then we define:
\[\rt(k,G)=\{x\in\cl(k):\ \exists K/k\text{ tame, }\gal(K/k)\cong G, \st(K/k)=x\}.\]

\begin{defn}\label{A'groups}
We define $A'$-groups inductively:
\begin{enumerate}
\item Finite abelian groups are $A'$-groups.
\item If $\G$ is an $A'$-group and $H$ is finite abelian of order prime to that of $\G$, then $H\rtimes_\mu \G$ is an $A'$-group, for any action $\mu$ of $\G$ on $H$.
\item If $\G_1$ and $\G_2$ are $A'$-groups, then $\G_1\times \G_2$ is an $A'$-group.
\end{enumerate}
\end{defn}

In the following proposition we find a relation between $A'$-groups and more classical kinds of groups.

\begin{prop}
Every $A'$-group is a solvable $A$-group (an $A$-group is a finite group with the property that all of its Sylow subgroups are abelian).
\end{prop}

\begin{proof}
Since abelian groups are obviously solvable $A$-groups, we have only to prove that the property of being a solvable $A$-group is preserved by constructions 2 and 3 in Definition \ref{A'groups}.

If $\G$, $\G_1$ and $\G_2$ are solvable and $H$ is abelian, then $H\rtimes_\mu \G$ and $\G_1\times \G_2$ are clearly solvable.

If $\G$ is an $A$-group and $H$ is abelian of order prime to that of $\G$, then for any prime $l$ dividing the order of $H$ an $l$-Sylow subgroup of $H\rtimes_\mu \G$ must be a subgroup of $H$ and thus must be abelian. If $l$ divides the order of $\G$ then an $l$-Sylow subgroup of $H\rtimes_\mu \G$ is isomorphic to one of $\G$ and thus it is abelian, by hypothesis. So $H\rtimes_\mu \G$ is an $A$-group.

If $\G_1$ and $\G_2$ are $A$-groups, then for any prime $l$, an $l$-Sylow subgroup of $\G_1 \times \G_2$ is a direct product of $l$-Sylow subgroups of $\G_1$ and $\G_2$ and hence it is abelian, and $\G_1 \times \G_2$ is an $A$-group.
\end{proof}

It is an open question whether the converse of the proposition is true or not.

The main result we are going to prove is that the realizable classes for a number field $k$ and an $A'$-group $G$ of odd order form a group; this covers many cases not previously known. Further we use the same techniques to determine $\rt(k,D_n)$ for any odd integer $n$.

In contrast with many other papers on the subject, we systematically use class field theory (instead of Kummer theory and cyclotomic descent).

This paper is a slightly shortened version of parts of the author's PhD thesis \cite{tesi}. For earlier results see \cite{BrucheSodaigui}, \cite{ByottGreitherSodaigui}, \cite{Carter}, \cite{CarterSodaigui_quaternionigeneralizzati}, \cite{Endo}, \cite{GodinSodaigui_A4}, \cite{GodinSodaigui_ottaedri}, \cite{Long2}, \cite{Long0},  \cite{Massy}, \cite{McCulloh},  \cite{Sodaigui1}, \cite{Sodaigui2} and \cite{Soverchia}.

\end{section}

\markboth{Acknowledgements}{Acknowledgements}
\begin{section}*{Acknowledgements}
\addcontentsline{toc}{section}{Acknowledgements}
I am very grateful to Professor Cornelius Greither and to Professor Roberto Dvornicich for their advice and for the patience they showed, assisting me in the writing of my PhD thesis with a lot of suggestions and corrections. I also wish to thank the Scuola Normale Superiore of Pisa, for its role in my mathematical education and for its support during the time I was working on my PhD thesis.
\end{section}

\begin{section}{Preliminary results}
\begin{subsection}{Steinitz classes}
We start recalling some well-known results about Steinitz classes.

It is a classical result (which can be deduced by Propositions 8 and 14 of chapter III of \cite{Lang}) that the discriminant of a tamely ramified Galois extension $K/k$ of number fields is
\[\disc(K/k)=\prod_{\p}\p^{(e_\p-1)\frac{[K:k]}{e_\p}}\]
where $e_\p$ is the ramification index of $\p$.

The discriminant is closely related to the Steinitz class by the following theorem.

\begin{teo}\label{stdisc}
Assume $K$ is a finite Galois extension of a number field $k$.
\begin{enumerate}
\item[(a)] If its Galois group either has odd order or has a noncyclic $2$-Sylow subgroup then $\disc(K/k)$ is the square of an ideal and this ideal represents the Steinitz class of the extension.
\item[(b)] If its Galois group is of even order with a cyclic $2$-Sylow subgroup and $\alpha$ is any element of $k$ whose square root generates the quadratic subextension of $K/k$ then $\disc(K/k)/\alpha$ is the square of a fractional ideal and this ideal represents the Steinitz class of the extension.
\end{enumerate}
\end{teo}

\begin{proof}
This is a corollary of Theorem I.1.1 in \cite{Endo}. In particular it is shown in \cite{Endo} that in case (b) $K/k$ does have exactly one quadratic subextension.
\end{proof}

Further, considering Steinitz classes in towers of extensions, we will need the following proposition.

\begin{prop}\label{stintermediateextension}
Suppose $K/E$ and $E/k$ are number fields extensions. Then
\[\st(K/k)=\st(E/k)^{[K:E]}\norm_{E/k}(\st(K/E)).\]
\end{prop}

\begin{proof}
This is Proposition I.1.2 in \cite{Endo}.
\end{proof}

\end{subsection}

\begin{subsection}{Class field theory}
To prove our results we will use techniques from class field theory. We will use the notations and some results of \cite{Neukirch}. 

We will denote by $U_\p$ the units of a number field $k$ completed at a prime $\p$, by $I_k$ the idele group of $k$ and by $C_k=I_k/k^*$ the idele class group. We set
\[U_\p^{n}=\begin{cases}U_\p&\text{if $n=0$}\\ 1+\p^{n} &\text{if $n>0$}.\end{cases}\]
For any cycle $\mathfrak m=\prod_\p \p^{n_\p}$ we consider the groups
$I_k^\mathfrak m=\prod_\p U_\p^{n_\p}$ and the congruence subgroup mod $\mathfrak m$ of $C_k$, i.e. $C_k^\mathfrak m=I_k^\mathfrak m\cdot k^*/k^*\subseteq C_k$.

For every prime $\p$ we have the canonical injection
\[[\ ]:k_\p^*\to C_k,\]
which associates to $a_\p\in k_\p^*$ the class of the idele
\[[a_\p]=(\dots,1,1,1,a_\p,1,1,1,\dots).\]

To construct number fields extensions with a given Steinitz class we will use the following results.

\begin{teo}\label{classfieldtheory}
Let $G$ be an abelian group. Every surjective homomorphism $\varphi:C_k\to G$ whose kernel contains a congruence subgroup $C_k^\mathfrak m$ is the norm residue symbol of a unique extension $K/k$ with Galois group isomorphic to $G$ and 
$\varphi([U_\p])$ is its inertia group for the prime $\p$. In particular
\[e_\p(K/k)=\#\varphi([U_\p])\]
and if the primes dividing the order of $G$ do not divide $\mathfrak m$, then the extension is tame.
\end{teo}

\begin{proof}
By Theorem IV.7.1 of \cite{Neukirch} there exists a unique abelian extension $K/k$ with $\norm_{K/k}C_K=\ker\varphi$. By Theorem IV.6.5 in \cite{Neukirch} the global residue symbol of $K/k$ gives an isomorphism $C_k/\ker\varphi=C_k/\norm_{K/k}C_K\to\gal(K/k)^\mathrm{ab}=\gal(K/k)$ and thus clearly $\gal(K/k)\cong G$. Now let $K_1$ and $K_2$ be two fields corresponding to the same residue symbol, then $\norm_{K_1/k}C_{K_1}=\norm_{K_2/k}C_{K_2}$ and so, by Theorem IV.7.1 of \cite{Neukirch}, $K_1=K_2$.

The group $\varphi([U_\p])$ is the inertia group for the prime $\p$ because of Theorem III.8.10 in \cite{Neukirch} and Proposition IV.6.6 in \cite{Neukirch}.
\end{proof}

Let $L/K$ and $K/k$ be an abelian and a Galois extension of number fields respectively, such that $L/k$ is normal, $\mathcal U=\gal(L/K)$ and $\Delta=\gal(K/k)$. Let $\delta\in\Delta$, $\sigma\in\mathcal U$ and let $\tilde\delta,\tilde \delta'\in\gal(L/k)$ be two extensions of $\delta$ to $\gal(L/k)$. Then $\tilde\delta'^{-1}\tilde\delta\in\mathcal U$ and, by the commutativity of $\mathcal U$, we have that
\[\tilde\delta_*\sigma=\tilde\delta\sigma\tilde\delta^{-1}=\tilde\delta'\tilde\delta'^{-1}\tilde\delta\sigma\tilde\delta^{-1}\tilde\delta'\tilde\delta'^{-1}=\tilde\delta'\sigma \tilde\delta'^{-1}=\tilde\delta'_*\sigma,\]
so that we can define $\delta_*:\mathcal U\to\mathcal U$ by $\delta_*=\tilde\delta_*$. 

\begin{prop}\label{Greither}
Let $K/k$ be a finite tame extension with Galois group $\Delta$, let $\mathcal U$ be a finite abelian group and let $\phi:\Delta\to\mathrm{Aut}(\mathcal U)$ be an action of $\Delta$ on $\mathcal U$. Then for a $\Delta$-invariant surjective homomorphism $\varphi:C_K\to\mathcal U$, whose kernel contains a congruence subgroup $C_K^\mathfrak m$, the extension $L/K$ given by Theorem \ref{classfieldtheory} is Galois over $k$. The following sequence is exact
\[1\to\mathcal U\to\gal(L/k)\to\Delta\to 1\]
and the induced action of $\Delta$ on $\mathcal U$ is the given one.
\end{prop}

\begin{proof}
Let $\tilde K$ be the maximal abelian extension of $K$; by standard arguments $\tilde K/k$ is Galois. Since $\tilde K\supset L$, there is a normal closure $L_1$ of $L/k$ in $\tilde K$ and the extension $L_1/K$ is finite and abelian. Let $\pi:\gal(L_1/K)\to\gal(L/K)$ be the projection, then $L$ is the fixed field of $\ker\pi$. By Proposition II.3.3 in \cite{Neukirch}
\[\pi=(\ ,L/K)\circ r_{L_1/K}=\varphi\circ r_{L_1/K}\]
and for $\delta\in\gal(L_1/k)$ we have, using also the hypothesis of $\Delta$-invariance,
\[\delta_*\circ\pi=\delta_*\circ\varphi\circ r_{L_1/K}=\varphi\circ\delta\circ r_{L_1/K}=\varphi\circ r_{L_1/K}\circ\delta_*=\pi\circ\delta_*.\]
Thus
\[\delta_*\ker\pi=\ker(\pi\circ\delta_*^{-1})=\ker(\delta_*^{-1}\circ\pi)=\ker\pi.\]
So $\ker\pi$ is normal in $\gal(L_1/k)$. It follows that $L/k$ is Galois. The exactness of the sequence is obvious and the statement about the action of $\Delta$ on $\mathcal U$ follows from Proposition II.3.3 in \cite{Neukirch}, since the given action is the only one for which the diagram on the right commutes.
\end{proof}

For any cycle $\mathfrak m=\prod_\p \p^{n_\p}$ we call $H_{K/k}^\mathfrak m=\norm_{K/k}J_K^\mathfrak m\cdot P_k^\mathfrak m$, where  $J_k^\mathfrak m$ is the group of all ideals prime to $\mathfrak m$ and $P_k^\mathfrak m$ is the group of all principal ideals generated by an element $a\equiv 1\pmod{\p^{n_\p}}$ for all $\p|\mathfrak m$.

Let $K$ be an abelian extension of $k$, contained in the ray class field mod $\mathfrak m$; the cycle $\mathfrak m$ is called a \emph{cycle of declaration} for $K/k$.

\begin{prop}\label{ideallattice}
Let $K,K_1,K_2$ be finite abelian extensions of a number field $k$ and let $\mathfrak m$ be a cycle of declaration for them. Then there is an isomorphism
\[\pi_\mathfrak m:\norm_{K/k}C_K/C_k^\mathfrak m\to H_{K/k}^\mathfrak m/P_k^\mathfrak m,\]
and
\[K_1\subseteq K_2\iff H_{K_1/k}^\mathfrak m\supseteq H_{K_2/k}^\mathfrak m,\]
\[H_{K_1\cdot K_2/k}^\mathfrak m=H_{K_1/k}^\mathfrak m\cap H_{K_2/k}^\mathfrak m,\quad H_{K_1\cap K_2/k}^\mathfrak m=H_{K_1/k}^\mathfrak m\cdot H_{K_2/k}^\mathfrak m.\]
\end{prop}

\begin{proof}
By Proposition IV.8.1 in \cite{Neukirch}, there exists a surjective homomorphism $\pi_\mathfrak m:C_k\to J_k^\mathfrak m/P_k^\mathfrak m$ and by the exact commutative diagram in Theorem IV.8.2 in \cite{Neukirch}, we obtain that $\pi_\mathfrak m(\norm_{K/k}C_K)=H_{K/k}^\mathfrak m/P_k^\mathfrak m$. By Theorem IV.7.1 of \cite{Neukirch}, $\norm_{K/k}C_K\supseteq C_k^\mathfrak m$ ($\mathfrak m$ is a cycle of declaration of $K/k$) and then by Proposition IV.8.1 in \cite{Neukirch} it is the kernel of $\pi_\mathfrak m:\norm_{K/k}C_K\to H_{K/k}^\mathfrak m/P_k^\mathfrak m$.

Now the result follows by Theorem IV.7.1 of \cite{Neukirch} and by the fact that $H_{K/k}^\mathfrak m$ is the counterimage of $H_{K/k}^\mathfrak m/P_k^\mathfrak m$ by the projection $J_k^\mathfrak m\to J_k^\mathfrak m/P_k^\mathfrak m$.
\end{proof}

\begin{prop}\label{Artinrayclassfield}
Let $k^\mathfrak m$ be the ray class field modulo a cycle $\mathfrak m$ of a number field $k$. Then
\[\left(\frac{k^\mathfrak m/k}{\cdot}\right):J_k^\mathfrak m/P_k^\mathfrak m\to\gal(k^\mathfrak m/k)\]
is an isomorphism.
\end{prop}

\begin{proof}
By definition of the ray class field mod $\mathfrak m$, $\norm_{k^\mathfrak m/k}C_{k^\mathfrak m}=C_{k}^\mathfrak m$ and thus by Proposition \ref{ideallattice} we obtain that $H_{k^\mathfrak m/k}^\mathfrak m/P_k^\mathfrak m$ is the trivial group. We conclude using Theorem IV.8.2 in \cite{Neukirch}.
\end{proof}

\begin{prop}\label{Tschebotarev2}
Let $\mathfrak m$ be a cycle for a number field $k$. Then each class in the ray class group modulo $\mathfrak m$ contains infinitely many prime ideals of absolute degree $1$.
\end{prop}

\begin{proof}
For each ray class in $J_k^\mathfrak m/P_k^\mathfrak m$ we can consider the automorphisms $\sigma\in\gal(k^\mathfrak m/k)$ corresponding to it by the isomorphism of Proposition \ref{Artinrayclassfield}. By Chebotarev Theorem (V.6.4 in \cite{Neukirch}), there exist infinitely many prime ideals $\p$ in $k$, unramified in $k^\mathfrak m$, of absolute degree $1$ and with $\sigma=\left(\frac{k^\mathfrak m/k}{\p}\right)$. By construction they must be in the given ray class.
\end{proof}

\begin{defn}\label{defW}
Let $K/k$ be a finite abelian extension of number fields and let $\mathfrak m$ be a cycle of declaration of $K/k$. We define
\[W(k,K)=\norm_{K/k}J_K^\mathfrak m\cdot P_k/P_k=H_{K/k}^\mathfrak m\cdot P_k/P_k.\]
If $\zeta_m$ is an $m$-th root of unity we use the notation $W(k,m)=W(k,k(\zeta_m))$.
\end{defn}

\begin{prop}\label{Wclassfield}
By class field theory $W(k,K)$ corresponds to the maximal unramified subextendion of $K/k$, i.e.
\[W(k,K)=H_{K\cap k^1/k}^1/P_k,\]
where $k^1$ is the Hilbert class field of $k$. In particular $W(k,K)$ does not depend on the choice of the cycle of declaration $\mathfrak m$ of $K/k$.
\end{prop}

\begin{proof}
By Theorem IV.8.2 of \cite{Neukirch} and by Proposition \ref{Artinrayclassfield} the kernel of \[\left(\frac{k^1/k}{\cdot}\right):J_k^\mathfrak m/P_k^\mathfrak m\to\gal(k^1/k)\]
is $H_{k^1/k}^\mathfrak m/P_k^\mathfrak m=(P_k\cap J_k^\mathfrak m)/P_k^\mathfrak m$, i.e. $H_{k^1/k}^\mathfrak m=P_k\cap J_k^\mathfrak m$
and, by Proposition \ref{ideallattice},
\[H_{K\cap k^1/k}^\mathfrak m=H_{K/k}^\mathfrak m\cdot H_{k^1/k}^\mathfrak m=H_{K/k}^\mathfrak m\cdot (P_k\cap J_k^\mathfrak m).\]

Let $x\in H_{K\cap k^1/k}^1/P_k$, then by Proposition \ref{Tschebotarev2} there exists a prime $\p\nmid\mathfrak m$ in the class of $x$, i.e., recalling also the definition of $H_{K\cap k^1/k}^\mathfrak m$,
\[\p\in H_{K\cap k^1/k}^1\cap J_k^\mathfrak m=H_{K\cap k^1/k}^\mathfrak m\cdot P_k=H_{K/k}^\mathfrak m\cdot (P_k\cap J_k^\mathfrak m)\cdot P_k=H_{K/k}^\mathfrak m\cdot P_k\]
and so $x\in H_{K/k}^\mathfrak m\cdot P_k/P_k$. Thus
\[H_{K\cap k^1/k}^1/P_k\subseteq H_{K/k}^\mathfrak m\cdot P_k/P_k=H_{K\cap k^1/k}^\mathfrak m/P_k\]
and the opposite inclusion is trivial.

Thus we have proved that
\[W(k,K)=H_{K/k}^\mathfrak m\cdot P_k/P_k=H_{K\cap k^1/k}^1/P_k.\]
\end{proof}

The following results are similar to the characterizations of $W(k,K)$ given in \cite{Endo}.
\begin{prop}\label{I.2.2Endo}
Let $K/k$ be a finite abelian extension of number fields. Then the following subsets of the class group of $k$ are equal to $W(k,K)$:
\[\begin{split}
&W_1=\{x\in J_k/P_k:\text{$x$ contains infinitely many primes of absolute degree $1$}\\&\qquad\ \ \text{ splitting completely in $K$}\}\\
&W_2=\{x\in J_k/P_k:\text{$x$ contains a prime splitting completely in $K$}\}\\
&W_3=\norm_{K/k}(J_K)\cdot P_k/P_k.\end{split}\]
\end{prop}

\begin{proof}
Let $x\in W(k,K)$ and let $\mathfrak m$ be a cycle of declaration of $K/k$. By definition $x=\mathfrak a\cdot P_k$, where $\mathfrak a\in  H_{K/k}^\mathfrak m$. By Proposition \ref{Tschebotarev2} there exist infinitely many primes of absolute degree $1$ in the ray class modulo $\mathfrak m$ containing $\mathfrak a$; let $\p$ be one of them, which does not ramify in $K/k$. Then $\p=\mathfrak a\cdot (b)$, where $(b)\in P_k^\mathfrak m$, and thus $\p\in H_{K/k}^\mathfrak m$ and by Theorem IV.8.4 in \cite{Neukirch} we can conclude that $\p$ splits completely in $K$. Thus $x\in W_1$ and we have proved that $W(k,K)\subseteq W_1$.

Obviously $W_1\subseteq W_2$.

Let $x\in W_2$ and $\p$ be a prime in $x$ which splits completely in $K$. Then for any prime divisor $\P$ of $\p$ in $K$, $\norm_{K/k}(\P)=\p$. Thus $x=\norm_{K/k}(\P)\cdot P_k$ and hence $W_2\subseteq W_3$.

Recalling Proposition \ref{Wclassfield} we obtain that
\[N_{K/k}(J_K)\cdot P_k/P_k\subseteq N_{K\cap k^1/k}(J_{K\cap k^1})\cdot P_k/P_k=H_{K\cap k^1/k}^1/P_k=W(k,K).\]
\end{proof}

In the case of cyclotomic extensions we obtain some further results.
\begin{lemma}\label{Wkpnorme}
Let $\p$ be a prime in $k$ of absolute degree $1$, splitting completely in $k(\zeta_m)$ and unramified over $\Q$. Then $\norm_{k/\Q}(\p)\in P_\Q^\mathfrak m$, where $\mathfrak m=m\cdot p_\infty$.
\end{lemma}

\begin{proof}
By hypothesis $\oo_k/\p$ is the finite field with $p$ elements, where $\norm_{k/\Q}(\p)=(p)$, and $(\oo_k/\p)^*$ contains a primitive $m$-th root of unity, i.e. an element of order $m$. Hence $m$ must divide $|\oo_k/\p|-1=p-1$, i.e. $p\equiv 1\pmod m$, which is equivalent to the assertion.
\end{proof}

\begin{lemma}\label{Wram2}
Let $k$ be a number field, let $\mathfrak m$ be a cycle of declaration of $k(\zeta_m)/k$ and let $\mathfrak a\in J_k^\mathfrak m$ be such that $\norm_{k/\Q}(\mathfrak a)\in P_\Q^\mathfrak m$, then $\mathfrak a\in H_{k(\zeta_m)/k}^\mathfrak m$, i.e. the class of $\mathfrak a$ is in $W(k,m)$.
\end{lemma}

\begin{proof}
By Proposition II.3.3 in \cite{Neukirch},
\[\left.\left(\frac{k(\zeta_m)/k}{\mathfrak a}\right)\right|_{\Q(\zeta_m)}=\left(\frac{\Q(\zeta_m)/\Q}{\norm_{k/\Q}(\mathfrak a)}\right)=1.\]
Of course also the restriction of $\left(\frac{k(\zeta_m)/k}{\mathfrak a}\right)$ to $k$ is trivial; thus we have proved that
\[\left(\frac{k(\zeta_m)/k}{\mathfrak a}\right)=1,\]
i.e. that $\mathfrak a\in H_{k(\zeta_m)/k}^\mathfrak m$.
\end{proof}

\begin{lemma}\label{Wram}
Let $K/k$ be a tamely ramified abelian extension of number fields and let $\p$ be a prime ideal in $k$ whose ramification index in $K/k$ is $e$, then $\norm_{k/\Q}(\p)\in P_\Q^\mathfrak m$, where $\mathfrak m=e\cdot p_\infty$. In particular, by Lemma \ref{Wram2}, $\p\in H_{k(\zeta_e)/k}^\mathfrak m$ and so its class is in $W(k,e)$.
\end{lemma}
\begin{proof}
This is Lemma I.2.1 of \cite{Endo}.
\end{proof}

\end{subsection}

\end{section}

\begin{section}{Main results}
Let $\G$ be a finite group of order $m$, let $H=C(n_1)\times\dots\times C(n_r)$ be an abelian group of order $n$, with generators $\tau_1,\dots,\tau_r$ and with $n_{i+1}|n_i$. Let
\[\mu:\G\to \aut(H)\]
be an action of $\G$ on $H$ and let
\[0\to H\stackrel\varphi\longrightarrow G\stackrel\psi\longrightarrow \G\to 0\]
be an exact sequence of groups such that the induced action of $\G$ on $H$ is $\mu$. We assume that the group $G$ is determined, up to isomorphism, by the above exact sequence and by the action $\mu$ (this is true e.g. if $|\G|$ and $|H|$ are coprime). We are going to study $\rt(k,G)$. 

We define
\[\eta_G=\begin{cases}1&\text{if some (and hence every) $2$-Sylow subgroup of $G$ is not cyclic}\\2&\text{if some (and hence every) $2$-Sylow subgroup of $G$ is cyclic}\end{cases}\]
and in a similar way we define $\eta_H$ and $\eta_\G$. We will always use the letter $l$ only for prime numbers, even if not explicitly indicated.

We say that $(K,k_1,k)$ is of type $\mu$ if $k_1/k$, $K/k_1$ and $K/k$ are Galois extensions with Galois groups isomorphic to $\G$, $H$ and $G$ respectively and such that the action of $\gal(k_1/k)\cong \G$ on $\gal(K/k_1)\cong H$ is given by $\mu$. For any $\G$-extension $k_1$ of $k$ we define $\rt(k_1,k,\mu)$ as the set of those ideal classes of $k_1$ which are Steinitz classes of a tamely ramified extension $K/k_1$ for which $(K,k_1,k)$ is of type $\mu$.

\begin{subsection}{Some definitions and simple properties}

For any $\tau\in H$ we define 
\[\tilde\mu_{k,\mu,\tau}:\G\times\gal(k(\zeta_{o(\tau)})/k)\to\aut(H)\]
by $\tilde\mu_{k,\mu,\tau}((g_1,g_2))=\mu(g_1)$ for any $(g_1,g_2)\in \G\times\gal(k(\zeta_{o(\tau)})/k)$ and
\[\tilde\nu_{k,\mu,\tau}:\G\times\gal(k(\zeta_{o(\tau)})/k)\to(\Z/o(\tau)\Z)^*\]
by $\tilde\nu_{k,\mu,\tau}((g_1,g_2))=\nu_{k,\tau}(g_2)$ where $g_2(\zeta_{o(\tau)})=\zeta_{o(\tau)}^{\nu_{k,\tau}(g_2)}$ for any $(g_1,g_2)\in \G\times\gal(k(\zeta_{o(\tau)})/k)$. Let
\[\begin{split}\tilde G_{k,\mu,\tau}&=\left\{g\in \G\times\gal(k(\zeta_{o(\tau)})/k):\ \tilde\mu_{k,\mu,\tau}(g)(\tau)=\tau^{\tilde\nu_{k,\mu,\tau}(g)}\right\}
\\&=\left\{(g_1,g_2)\in \G\times\gal(k(\zeta_{o(\tau)})/k):\ \mu(g_1)(\tau)=\tau^{\nu_{k,\tau}(g_2)}\right\}.\end{split}\]

We define
\[G_{k,\mu,\tau}=\left\{g\in\gal(k(\zeta_{o(\tau)})/k):\ \exists g_1\in\G,\ (g_1,g)\in \tilde G_{k,\mu,\tau}\right\}\]
and $E_{k,\mu,\tau}$ as the fixed field of $G_{k,\mu,\tau}$ in $k(\zeta_{o(\tau)})$.

\begin{lemma}\label{Gtilde}
For any $\tau\in H$, $G_{k,\mu,\tau}$ is a subgroup of $\gal(k(\zeta_{o(\tau)})/k)$.
\end{lemma}
\begin{proof}

If $(g_1,g_2),(\tilde g_1,\tilde g_2)\in \tilde G_{k,\mu,\tau}$, then 
\[\begin{split}\tau^{\tilde\nu_{k,\mu,\tau}((g_1\tilde g_1,g_2\tilde g_2))}&=\tau^{\nu_{k,\tau}(g_2)\nu_{k,\tau}(\tilde g_2)}=\mu(g_1)\left(\tau^{\nu_{k,\tau}(\tilde g_2)}\right)\\&=\mu(g_1)(\mu(\tilde g_1)(\tau))=\tilde\mu_{k,\mu,\tau}((g_1\tilde g_1,g_2\tilde g_2))(\tau)\end{split}\]
and
\[\begin{split}\tau^{\tilde\nu_{k,\mu,\tau}\left(\left(g_1^{-1},g_2^{-1}\right)\right)}&=\tau^{\nu_{k,\tau}(g_2)^{-1}}=\mu\left(g_1^{-1}\right)\left(\mu(g_1)\left(\tau^{\nu_{k,\tau}(g_2)^{-1}}\right)\right)\\&=\tilde\mu_{k,\mu,\tau}\left(\left(g_1^{-1},g_2^{-1}\right)\right)(\tau).\end{split}\]

Hence $(g_1\tilde g_1,g_2\tilde g_2)$,$\left(g_1^{-1},g_2^{-1}\right)\in \tilde G_{k,\mu,\tau}$ and the set $G_{k,\mu,\tau}$ is a subgroup of $\gal(k(\zeta_{o(\tau)})/k)$.
\end{proof}

Given a $\G$-extension $k_1$ of $k$, there is an injection of $\gal(k_1(\zeta_{o(\tau)})/k)$ into $\G\times\gal(k(\zeta_{o(\tau)})/k)$ (defined in the obvious way). We will always identify $\gal(k_1(\zeta_{o(\tau)})/k)$ with its image in $\G\times\gal(k(\zeta_{o(\tau)})/k)$. So we may consider the subgroup
\[\tilde G_{k_1/k,\mu,\tau}=\tilde G_{k,\mu,\tau}\cap \gal(k_1(\zeta_{o(\tau)})/k)\]
of $\tilde G_{k,\mu,\tau}$. 
Let $Z_{k_1/k,\mu,\tau}$ be its fixed field in $k_1(\zeta_{o(\tau)})$.

If $k_1\cap k(\zeta_{o(\tau)})=k$ then $\gal(k_1(\zeta_{o(\tau)})/k)\cong \G\times\gal(k(\zeta_{o(\tau)})/k)$ and hence 
$\tilde G_{k_1/k,\mu,\tau}=\tilde G_{k,\mu,\tau}$.

\begin{lemma}\label{defZlgen}
For any $\tau\in H$, $k_1Z_{k_1/k,\mu,\tau}=k_1(\zeta_{o(\tau)})$.
\end{lemma}

\begin{proof}
Let $g\in \gal(k_1(\zeta_{o(\tau)})/k_1)\cap \tilde G_{k_1/k,\mu,\tau}$,
then $g|_{k_1}=1$, i.e. $\tilde\mu_{k,\mu,\tau}(g)(\tau)=\tau$, and $\tau^{\tilde\nu_{k,\mu,\tau}(g)}=\tilde\mu_{k,\mu,\tau}(g)(\tau)=\tau$. Thus $g(\zeta_{o(\tau)})=\zeta_{o(\tau)}$ and we conclude that $g=1$. We have proved that
\[\gal(k_1(\zeta_{o(\tau)})/k_1)\cap \tilde G_{k_1/k,\mu,\tau}=1\]
i.e. that 
\[k_1Z_{k_1/k,\mu,\tau}=k_1(\zeta_{o(\tau)}).\]
\end{proof}

\begin{lemma}\label{caratterizzaE}
Let $\tau\in H$, then 
\[E_{k,\mu,\tau}\subseteq Z_{k_1/k,\mu,\tau}\cap k(\zeta_{o(\tau)})\]
and we have an equality if $k_1\cap k(\zeta_{o(\tau)})=k$.
\end{lemma}

\begin{proof}
We observe that
\[\begin{split}
G_{k,\mu,\tau}&\supseteq \left\{g_2\in\gal(k(\zeta_{o(\tau)})/k):\ \exists g_1\in\G,\ (g_1,g_2)\in \tilde G_{k_1/k,\mu,\tau}\right\}\\
&=\mathrm{res}^{k_1(\zeta_{o(\tau)})}_{k(\zeta_{o(\tau)})}(\tilde G_{k_1/k,\mu,\tau})\\
&=\mathrm{res}^{k_1(\zeta_{o(\tau)})}_{k(\zeta_{o(\tau)})}(\tilde G_{k_1/k,\mu,\tau})\mathrm{res}^{k_1(\zeta_{o(\tau)})}_{k(\zeta_{o(\tau)})}(\gal(k_1(\zeta_{o(\tau)})/k(\zeta_{o(\tau)})))\\
&=\mathrm{res}^{k_1(\zeta_{o(\tau)})}_{k(\zeta_{o(\tau)})}(\gal(k_1(\zeta_{o(\tau)})/Z_{k_1/k,\mu,\tau}\cap k(\zeta_{o(\tau)})))\\
&=\gal(k(\zeta_{o(\tau)})/Z_{k_1/k,\mu,\tau}\cap k(\zeta_{o(\tau)}))
\end{split}\]

i.e. that
\[E_{k,\mu,\tau}\subseteq Z_{k_1/k,\mu,\tau}\cap k(\zeta_{o(\tau)}).\]
If $k_1\cap k(\zeta_{o(\tau)})=k$ then $\tilde G_{k_1/k,\mu,\tau}=\tilde G_{k,\mu,\tau}$ and we have equalities.
\end{proof}

\begin{lemma}\label{inclusioniWEWZ}
Let $\tau\in H$, then 
\[W(k,Z_{k_1/k,\mu,\tau})\subseteq W(k,E_{k,\mu,\tau}).\]
If $k_1\cap k(\zeta_{o(\tau(l))})=k$ and every subextension of $k_1/k$ is ramified then
\[W(k,Z_{k_1/k,\mu,\tau})=W(k,E_{k,\mu,\tau}).\]
\end{lemma}

\begin{proof}
By Lemma \ref{caratterizzaE} it is obvious that \[W(k,Z_{k_1/k,\mu,\tau})\subseteq W(k,E_{k,\mu,\tau}).\]
Now we assume that $k_1/k$ has no unramified subextensions and we prove that
\[k^1\cap k_1(\zeta_{o(\tau)})\subseteq k(\zeta_{o(\tau)}),\]
where $k^1$ is the ray class field modulo $1$, i.e. the Hilbert class field.
If that is not true, then $k(\zeta_{o(\tau)})\subsetneq(k^1\cap k_1(\zeta_{o(\tau)}))\cdot k(\zeta_{o(\tau)})\subseteq k_1(\zeta_{o(\tau)})$ and the extension $(k^1\cap k_1(\zeta_{o(\tau)}))\cdot k(\zeta_{o(\tau)})/k(\zeta_{o(\tau)})$ is ramified at a prime ramified in $k_1/k$. This prime must ramify also in $k^1\cap k_1(\zeta_{o(\tau)})/k$, which is impossible. Thus if $k_1\cap k(\zeta_{o(\tau)})=k$ and $k_1/k$ has no unramified subextensions then, recalling also Lemma \ref{caratterizzaE},
\[k^1\cap E_{k,\mu,\tau}=k^1\cap Z_{k_1/k,\mu,\tau}\cap k(\zeta_{o(\tau)})=k^1\cap Z_{k_1/k,\mu,\tau}\cap k_1(\zeta_{o(\tau)})=k^1\cap Z_{k_1/k,\mu,\tau}\]
and by Proposition \ref{Wclassfield} we conclude that $W(k,E_{k,\mu,\tau})= W(k,Z_{k_1/k,\mu,\tau})$.
\end{proof}

\end{subsection}
\begin{subsection}{Some realizable classes for nonabelian groups}

First of all we need a more general version of the Multiplication Lemma on page 22 in \cite{Endo} by Lawrence P. Endo.

\begin{lemma}\label{multiplicationmetacyclic}
Let $(K_1,k_1,k)$ and $(K_2,k_1,k)$ be extensions of type $\mu$, such that $(\disc(K_1/k_1),\disc(K_2/k_1))=1$ and $K_1/k_1$ and $K_2/k_1$ have no non-trivial unramified subextensions. Then there exists an extension $(K,k_1,k)$ of type $\mu$, such that $K\subseteq K_1K_2$ and for which
\[\st(K/k_1)=\st(K_1/k_1)\st(K_2/k_1).\]
\end{lemma}

\begin{proof}
The hypotheses of the lemma imply that $K_1$ and $K_2$ are linearly disjoint over $k_1$.
Let us fix isomorphisms such that the action of $\G\cong \gal(k_1/k)$ on $H\cong\gal(K_i/k_1)$ given by conjugation coincides with $\mu$. Let us embed $H$ into $\gal(K_1K_2/k_1)$ by means of the corresponding diagonal map
\[\mathrm{diag}:H\to\gal(K_1/k_1)\times\gal(K_2/k_1)\cong\gal(K_1K_2/k_1).\]
Let $K$ be the fixed field of $\mathrm{diag}(H)$. Then, by Endo's Multiplication Lemma (page 22 in \cite{Endo}), we know that $\gal(K/k_1)\cong H$ and that
\[\st(K/k_1)=\st(K_1/k_1)\st(K_2/k_1).\]

The action of $\G\cong \gal(k_1/k)$ on \[\gal(K_1K_2/k_1)\cong\gal(K_1/k_1)\times\gal(K_2/k_1)\]
is given by
\[\tilde\mu(g)((h_1,h_2))=(\mu(g)(h_1),\mu(g)(h_2)).\]
It follows that the action of $\G\cong \gal(k_1/k)$ on
\[\gal(K/k_1)=\gal(K_1K_2/k_1)/\mathrm{diag}(H)\cong H\]
(where the last isomorphism is given by the projection on the first component) coincides with the action $\mu$. Hence $(K,k_1,k)$ is of type $\mu$.
\end{proof}

For any integer $n\in\N$ and any prime $l$, we denote by $n(l)$ the power of $l$ such that $n(l)|n$ and $l\nmid n/n(l)$. For any $\tau\in H$ and for any prime $l$ dividing the order $o(\tau)$ of $\tau$ we define the element
\[\tau(l)=\tau^\frac{o(\tau)}{o(\tau)(l)}\]
in the $l$-Sylow subgroup $H(l)$ of $H$. From now on we will assume that $H$ is of odd order.

We recall some definitions and a classical result.

\begin{defn}Let $R$ be a commutative ring, $G$ a finite group and $H$ a subgroup of $G$. The operation of \emph{restriction of scalars} from $R[G]$ to $R[H]$ assigns to each left $R[G]$-module $M$ a left $R[H]$-module $\res^G_H(M)$, whose underlying abelian group is still $M$ and such that for $h\in H$ and $m\in M$, $hm$ is obtained considering $h$ as an element of $G$.
\end{defn}

\begin{defn}
Let $R$ be a commutative ring, $G$ a finite group and $H$ a subgroup of $G$. The operation of \emph{induction} from $R[H]$-modules to $R[G]$-modules assigns to each left $R[H]$-module $L$ a left $R[G]$-module $\ind^G_H(L)$, given by
\[\ind^G_H (L)=R[G]\otimes_{R[H]}L.\]
\end{defn}

\begin{teo}[Frobenius reciprocity]\label{rappresentazioneindotta}
Let $H$ be a subgroup of a group $G$ and let $L$ be a left $R[H]$-module and $M$ a left $R[G]$-module. Then there exists an isomorphism of $R$-modules
\[\tau:\hom_{R[H]}(L,\res^G_H(M))\to\hom_{R[G]}(\ind^G_H(L),M).\]
This isomorphism is such that
\[(\tau f)(g\otimes l)=g\cdot f(l)\]
\end{teo}

\begin{proof}
This is Theorem 10.8 in \cite{CurtisReiner}. The explicit description of $\tau$ may be deduced from the proof.
\end{proof}

We will only use the above result with $R=\Z$.

Let $k_1/k$ be an extension of number fields with Galois group $\G$. Let $\P_1,\dots,\P_t$ be prime ideals in $\oo_{k_1}$, unramified over $p_1,\dots,p_t\in\N$, so that the classes $x_i$ of the $\P_i$ are generators of $\cl(k_1)$ (they exist because of Proposition \ref{Tschebotarev2}) and let $\P_i^{h_i}=(\alpha_i)$, where $h_i$ is the order of $x_i$.

We define the homomorphism (the so-called content map)
\[\pi:I_k\to J_k,\qquad \alpha\mapsto\prod_{\p\nmid\infty}\p^{v_\p(\alpha_\p)}.\]
Let $\pi_{\P_i}$ be a prime element in $(k_1)_{\P_i}$ and $y_i=[\pi_{\P_i}]\in I_{k_1}$, then $\pi(y_i)=\P_i$ and
\[ a_i=\frac{1}{\alpha_i}y_i^{h_i}\in \prod_\P U_\P\]
is congruent to $y_i^{h_i}$ mod $k_1^*$. 

For any $\delta\in\G$ let $b_{\delta,i}\in\prod_\P U_\P$ and $\lambda_{\delta,i,j}\in\Z$ (they exist thanks to the exactness of the sequence $1\to\prod_\P U_\P/U_{k_1}\to C_{k_1}\to \cl(k_1)\to 1$) be such that
\[\delta(y_i)=b_{\delta,i}\prod_{j=1}^t y_j^{\lambda_{\delta,i,j}}.\]

Let $\{u_1,\dots,u_{T}\}$ be the union of a system of generators of the abelian group $U_{k_1}$ with $\{a_1,\dots,a_t\}$ and $\bigcup_{\delta\in\G}\{b_{\delta,1},\dots,b_{\delta,t}\}$.

Let $\iota$ be the map from the class group of $k$ to the class group of $k_1$ induced by the map which pushes up ideals of $k$ to ideals of $k_1$. 
\begin{lemma}\label{Erweiterungslemma}
A group homomorphism $\varphi_0:(\prod_\p U_\p)/U_k\to G$ can be extended to $\varphi:C_k\to G$ if and only if for $j=1,\dots,t$, $\varphi_0(a_j)=g_j^{h_j}$ with $g_j\in G$. We can request also that $\varphi(y_j)=g_j$.
\end{lemma}

\begin{proof}
($\Rightarrow$) We have 
\[\varphi_0(a_j)=\varphi(y_j^{h_j})=\varphi(y_j)^{h_j}\in G^{h_j}.\]

($\Leftarrow$) Let us define
\[B_{k}=\left(\left(\prod_\p U_\p\right)/U_k\times\langle e_1,\dots,e_t\rangle\right)/\{e_j^{h_j}/a_j|j=1,\dots,t\}\]
where the second component in the direct product is a free abelian group. We may extend the inclusion $i:(\prod_\p) U_\p/U_k\hookrightarrow C_k$ to $B_{k}$ by $e_j\mapsto y_j$ and thus also the map $\pi\circ i:(\prod_\p U_\p)/U_k\to\cl(k)$ by $e_j\mapsto x_j$. We obtain the following commutative diagram
\[
\xymatrix{
1\ar[rr]^{}&&\left(\prod_\p U_\p\right)/U_k\ar[rr]^{}\ar[d]^{\mathrm{id}}& &B_{k}\ar[rr]^{}\ar@{->}[d]&&\cl(k)\ar[d]^{\mathrm{id}}\ar[rr]^{}&&1\\
1\ar[rr]^{}&&\left(\prod_\p U_\p\right)/U_k\ar[rr]^{}& &C_k\ar[rr]^{}&&\cl(k)\ar[rr]^{}&&1
}
\]
where the horizontal sequences are exact. It follows that $B_{k}\cong C_k$. Now we define $\tilde\varphi:B_{k}\to G$ by $\tilde\varphi(a)=\varphi_0(a)$ for $a\in \prod U_\p/U_k$ and $\tilde\varphi(e_j)=g_j$. This is a good definition since
\[\tilde\varphi\left(\frac{e_j^{h_j}}{a_j}\right)=\frac{g_j^{h_j}}{\varphi_0(a_j)}=1.\]
By the isomorphism between $B_{k}$ and $C_k$ we obtain the requested $\varphi:C_k\to G$. Since the restriction of the isomorphism $B_k\cong C_k$ to $\left(\prod_\p U_\p\right)/U_k$ is the identity map, it is clear that $\varphi$ is an extension of $\varphi_0$.
\end{proof}

The following lemma is a crucial technical result.

\begin{lemma}\label{metaconstructivelemma1}
Let $k_1$ be a tame $\G$-extension of $k$ and let $x\in W(k,k_1(\zeta_{n_1}))$. Then there exist tame extensions of $k_1$ of type $\mu$, whose Steinitz classes (over $k_1$) are $\iota(x)^{\alpha}$, where:
\[\alpha=\sum_{j=1}^{r} \frac{n_j-1}{2}\frac{n}{n_j}+\frac{n_1-1}{2}\frac{n}{n_1}.\]
In particular there exist tame extensions of $k_1$ of type $\mu$ with trivial Steinitz class. 

We can choose these extensions so that they are unramified at all infinite primes, that the discriminants are prime to a given ideal $I$ of $\oo_k$ and that all their proper subextensions are ramified.
\end{lemma}

\begin{proof}
By Proposition \ref{I.2.2Endo}, $x$ contains infinitely many primes $\q$ of absolute degree $1$ splitting completely in $k_1(\zeta_{n_1})$. Let $\q$ be any such prime and let $\q \oo_{k_1}=\prod_{\delta\in\G}\delta(\qq)$ be its decomposition in $k_1$, let $g_{\qq}$ be a generator of $\kappa_{\qq}^*=U_{\qq}/U_{\qq}^1$. 
Now $\delta$ gives an isomorphism from $\kappa_{\qq}^*$ to $\kappa_{\delta(\qq)}^*$ and so we may define a generator
\[g_{\delta(\qq)}=\delta(g_{\qq})\]
of $\kappa_{\delta(\qq)}^*$ for any $\delta\in\G$. We also define generators $g_\P$ of $\kappa_\P^*$ for all the other prime ideals and for any $a\in \prod_\P U_\P$ we define $\tilde h_{\P,a}\in\Z$, through $g_\P^{\tilde h_{\P,a}}\equiv a_\P\pmod{\P}$. 

For any prime $\delta(\qq)$, dividing a prime $\q$ of absolute degree $1$ splitting completely in $k_1(\zeta_{n_1})$, let $h_{\delta(\qq),a}$ be the class of $\tilde h_{\delta(\qq),a}$ modulo $n_1$ (since $\delta(\qq)$ is of absolute degree $1$, it follows by Lemma \ref{Wkpnorme} that the order of $g_{\delta(\qq)}$ is a multiple of $n_1$, i.e. that $h_{\delta(\qq),a}$ is well defined). The set of all the possible $m T$-tuples \[(h_{\delta(\qq),u_j})_{\delta\in\G;\ j=1,\dots,T}\]
is finite. Then it follows from the pigeonhole principle that there are infinitely many $\q$ corresponding to the same $m T$-tuple.

Let $\q_1,\dots,\q_{r+1}$ be $r+1$ such prime ideals and $\qq_1,\dots,\qq_{r+1}$ primes of $k_1$ dividing them. We can assume that they are distinct and that they are prime to a fixed ideal $I$ and to $\P_1,\dots,\P_t$.

Now let us define $\varphi_i:\kappa_{\qq_{i}}^*\to H$, posing
\[\varphi_i(g_{\qq_{i}})=\tau_i,\]
for $i=1,\dots,r$, and $\varphi_{r+1}:\kappa_{\qq_{r+1}}^*\to H$, posing
\[\varphi_{r+1}(g_{\qq_{r+1}})=(\tau_1\dots \tau_r)^{-1}.\]
Then we extend $\varphi_i$ to
\[\tilde\varphi_i:\ind^\G_{\langle 1\rangle}\kappa_{\qq_{i}}^*\cong\prod_{\delta\in\G}\kappa_{\delta(\qq_{i})}^*\to H\]
using Theorem \ref{rappresentazioneindotta}.

Now let us define $\varphi_0:\prod_{\qq}\kappa_{\qq}^*\to H$, posing 
\[\begin{cases}
\varphi_0|_{\kappa_{\delta(\qq_{i})}^*}=\tilde\varphi_i &\text{for $i=1,\dots,r+1$ and $\delta\in\G$}\\
\varphi_0|_{\kappa_{\P}^*}=1 &\text{for $\P\nmid \q_1,\dots,\q_{r+1}$}.\end{cases}\]
By construction $\varphi_0$ is $\G$-invariant and hence, for any $\delta\in\G$,
\[\varphi_0\left(\prod_{i=1}^{r+1}g_{\delta(\qq_{i})}\right)=\varphi_0\left(\delta\left(\prod_{i=1}^{r+1} g_{\qq_{i}}\right)\right)=\delta_*\varphi_0\left(\prod_{i=1}^{r+1} g_{\qq_{i}}\right)=\delta_*(1)=1.\]
It follows that $\varphi_0(u_j)=1$ for $j=1,\dots,T$ and thus in particular $\varphi_0$ is trivial on $U_{k_1}$, on the $a_1,\dots,a_t$ and on $b_{\delta,1},\dots,b_{\delta,t}$ for any $\delta\in\G$. This means that $\varphi_0$ is well defined on $(\prod_\P U_\P)/U_{k_1}$ and that $\varphi_0(a_j)=1$ and $\varphi_0(b_{\delta,j})=1$. Then it follows from Lemma \ref{Erweiterungslemma} that $\varphi_0$ can be extended to $\varphi:C_{k_1}\to G$; the kernel of $\varphi_0$ contains $I_{k_1}^\mathfrak m$, where $\mathfrak m=\prod_{j=1}^{r+1}\qq_j$, and so $C_k^\mathfrak m\subseteq \ker\varphi$.
We can also assume that $\varphi(y_j)=1$, for all $j$. It follows from Theorem \ref{classfieldtheory} that there is an $H$-Galois extension of $k_1$, ramifying only in the primes above $\q_1,\dots,\q_{r+1}$, with indices $n_j$ for $j\in\{1,\dots,r\}$ and $n_1$ for $j=r+1$.

Further the action of an element of $\G$ on one of the $y_j$ gives a combination of some $b_{\delta,i}$ and $y_j$, on which $\varphi$ is trivial. Recalling that $\varphi_0$ is $\G$-invariant, it follows that also the homomorphism $\varphi$ is $\G$-invariant and so by Proposition \ref{Greither} and by the fact that $G$ is identified by the exact sequence $1\to H\to G\to \G\to 1$ and by the action $\mu$ (Theorem 7.41 in \cite{Rotman}) we obtain an extension of type $\mu$.
Its discriminant is
\[\disc=\left(\prod_{i=1}^{r}\q_i^{(n_i-1)\frac{n}{n_i}}\right)\q_{r+1}^{(n_1-1)\frac{n}{n_1}}\oo_{k_1}\]
Since the order of $H$ is odd, by Theorem \ref{stdisc} the Steinitz class is $\iota(x)^{\alpha}$. It is immediate to verify the additional conditions.

\end{proof}

\begin{lemma}\label{metaconstructivelemma}
Let $k_1$ be a $\G$-extension of $k$, let $l$ be a prime dividing $n$, $\tau\in H(l)\setminus\{1\}$ and let $x$ be any class in $W(k,Z_{k_1/k,\mu,\tau})$. Then there exist extensions of $k_1$ of type $\mu$, whose Steinitz classes (over $k_1$) are $\iota(x)^{\alpha_{l,j}}$, where:
\[\alpha_{l,1}=(l-1)\frac{n}{l}; \tag{a}\]
\[\alpha_{l,2}=(o(\tau)-1)\frac{n}{o(\tau)}; \tag{b}\]
\[\alpha_{l,3}=\frac{3(l-1)}{2}.\tag{c}\]
We can choose these extensions so that they satisfy the additional conditions of Lemma \ref{metaconstructivelemma1}.
\end{lemma}

\begin{proof}
By Lemma \ref{metaconstructivelemma1} there exists an extension $K$ of $k_1$ of type $\mu$ with trivial Steinitz class and such that $K/k_1$ is unramified at all infinite primes, that its discriminant is prime to a given ideal $I$ of $\oo_k$ and that all its subextensions are ramified.

By Proposition \ref{I.2.2Endo}, $x$ contains infinitely many primes $\q$ of absolute degree $1$ splitting completely in $Z_{k_1/k,\mu,\tau}$. Those primes obviously split completely also in the extension $k_1(\zeta_{o(\tau)})=k_1Z_{k_1/k,\mu,\tau}$ (the equality holds by Lemma \ref{defZlgen}) of $k_1$. We can assume that they do not ramify in $k_1/k$, that they are prime to $l$ and, by the pigeonhole principle, that there are prime ideals $\qq$ in $k_1$, dividing the $\q$, and with a fixed decomposition group $D$, of order $f$, in $k_1/k$; let $\rho=m/f$. We choose a set $\Delta$ of representatives of the cosets $\delta D$, with $\delta\in\G$. Then $\q \oo_{k_1}=\prod_{\delta\in\Delta}\delta(\qq)$ are the decompositions of the primes $\q$ in $k_1$.

Let $g_{\qq}$ be a generator of $\kappa_{\qq}^*=U_{\qq}/U_{\qq}^1$. 
Now $\delta\in\Delta$ gives an isomorphism from $\kappa_{\qq}^*$ to $\kappa_{\delta(\qq)}^*$ and so we may define a generator
\[g_{\delta(\qq)}=\delta(g_{\qq})\]
of $\kappa_{\delta(\qq)}^*$ for any $\delta\in\Delta$.
We know that any $\delta\in D$ defines an automorphism of $\kappa_{\qq}^*$, of the form
\[\delta(g_{\qq})=g_{\qq}^{\lambda_{\qq,\delta}},\]
where $\lambda_{\qq,\delta}$ is an integer.
We can extend $\delta\in D$ to a $\tilde\delta\in \gal(k_1(\zeta_{o(\tau)})/k)$ in a way such that $\tilde\delta(\tilde\qq)=\tilde\qq$, where $\tilde\qq$ is a prime in $k_1(\zeta_{o(\tau)})$ above $\qq$ (it is enough to extend $\delta$ in some way and then to multiply it by an appropriate element of $\gal(k_1(\zeta_{o(\tau)})/k_1)$).
This element acts as a $\lambda_{\qq,\delta}$-th power on $\kappa_{\tilde\qq}^*=\kappa_{\qq}^*$ (the equality holds because $\qq$ splits completely in $k_1(\zeta_{o(\tau)})$). Thus, for $\delta\in D$, \[\zeta_{o(\tau)}^{\tilde\nu_{k,\mu,\tau}(\tilde\delta)}=\tilde\delta\left(\zeta_{o(\tau)}\right)\equiv \zeta_{o(\tau)}^{\lambda_{\qq,\delta}} \pmod{\tilde\qq}\]
and, recalling that the powers of $\zeta_{o(\tau)}$ are distinct modulo $\tilde\qq$ (since $\tilde\qq$ is prime to $l$ and thus to $o(\tau)$), 
\[\lambda_{\qq,\delta}\equiv \tilde\nu_{k,\mu,\tau}(\tilde\delta)\pmod{o(\tau)}.\]
Since the prime $\q$ splits completely in $Z_{k_1/k,\mu,\tau}$ and $\tilde\delta(\tilde\qq)=\tilde\qq$, we obtain that $\tilde\delta\in\gal(k_1(\zeta_{o(\tau)})/Z_{k_1/k,\mu,\tau})$ and hence
\[\mu(\delta)(\tau)=\tilde\mu_{k,\mu,\tau}(\tilde\delta)(\tau)=\tau^{\tilde\nu_{k,\mu,\tau}(\tilde\delta)}=\tau^{\lambda_{\qq,\delta}}.\]

Defining the $h_{\delta(\qq),u_j}$ as in the previous lemma, the set of all the possible $\rho T$-tuples
\[(h_{\delta(\qq),u_j})_{\delta\in\Delta;\ j=1,\dots,T}\]
is finite. Then it follows from the pigeonhole principle that there are infinitely many $\q$ corresponding to the same $\rho T$-tuple.

Let $\q_1,\q_2,\q_3$ be $3$ such prime ideals and let $\qq_1,\qq_2,\qq_3$ be primes of $k_1$ dividing them. We can assume that they are distinct, that they are prime to a fixed ideal $I$, to $\P_1,\dots,\P_t$ and to $\disc(K/k_1)$ and that they satisfy all the above requests.

\begin{enumerate}

\item[(a)] Now let us define $\varphi_i:\kappa_{\qq_{i}}^*\to H$, for $i=1,2$, posing 
\[\varphi_1(g_{\qq_1})=\tau^\frac{o(\tau)}{l}\]
and
\[\varphi_2(g_{\qq_2})=\tau^{-\frac{o(\tau)}{l}}.\]
For $\delta\in D$, we have
\[\mu(\delta)(\varphi_1(g_{\qq_1}))=\mu(\delta)\left(\tau^\frac{o(\tau)}{l}\right)=\tau^{\lambda_{\qq_1,\delta}\frac{o(\tau)}{l}}=\varphi_1(g_{\qq_1}^{\lambda_{\qq_1,\delta}})=\varphi_1(\delta(g_{\qq_1})).\]
Thus $\varphi_1$ is a $D$-invariant homomorphism and the same is true for $\varphi_2$.

Then, for $i=1,2$, we extend $\varphi_i$ to
\[\tilde\varphi_i:\ind^{\G}_{D}\kappa_{\qq_{i}}^*\cong\prod_{\delta\in\Delta}\kappa_{\delta(\qq_{i})}^*\to H\]
using Theorem \ref{rappresentazioneindotta} and we define $\varphi_0:\prod_{\P}\kappa_{\P}^*\to H$, posing  \[\begin{cases}
\varphi_0|_{\kappa_{\delta(\qq_{i})}^*}=\tilde\varphi_i &\text{for $i=1,2$  and $\delta\in\Delta$}\\
\varphi_0|_{\kappa_{\P}^*}=1 &\text{for $\P\nmid \q_1,\q_2$}.
\end{cases}\]

As in Lemma \ref{metaconstructivelemma1} we can extend $\varphi_0$ to a $\G$-invariant surjective homomorphism $\varphi:C_{k_1}\to H$, whose kernel contains a congruence subgroup of $C_{k_1}$ and hence this is true also for
\[\varphi\cdot (\ ,K/k_1):C_{k_1}\to H.\]
We can conclude that there exists an extension of type $\mu$, with discriminant
\[\disc(K/k_1)\left((\q_1\q_2)^{(l-1)\frac{n}{l}}\oo_{k_1}\right)\]
and so its Steinitz class is $\iota(x)^{\alpha_{l,1}}$.

\item[(b)] Now let us define $\varphi_i:\kappa_{\qq_{i}}^*\to H$, for $i=1,2$, posing 
\[\varphi_1(g_{\qq_1})=\tau\]
and
\[\varphi_2(g_{\qq_2})=\tau^{-1}\]
Exactly as in the first case we obtain an extension of type $\mu$ with discriminant
\[\disc(K/k_1)\left((\q_1\q_2)^{(o(\tau)-1)\frac{n}{o(\tau)}}\oo_{k_1}\right)\]
and so its Steinitz class is $\iota(x)^{\alpha_{l,2}}$.

\item[(c)] We define $\varphi_i:\kappa_{\qq_{i}}^*\to H$, for $i=1,2,3$, posing 
\[\varphi_1(g_{\qq_1})=\tau^\frac{o(\tau)}{l},\]
\[\varphi_2(g_{\qq_2})=\tau^\frac{o(\tau)}{l},\]
and
\[\varphi_3(g_{\qq_3})=\tau^{-\frac{2o(\tau)}{l}}\]
Now we obtain an extension of type $\mu$ with discriminant
\[\disc(K/k_1)\left((\q_1\q_2\q_3)^{(l-1)\frac{n}{l}}\oo_{k_1}\right)\]
and so its Steinitz class is $\iota(x)^{\alpha_{l,3}}$.
\end{enumerate}
Lemma \ref{metaconstructivelemma} is now completely proved.
\end{proof}

At this point we can prove the following proposition.

\begin{prop}\label{metaconstructive}
Let $l$ be a prime dividing the order $n$ of $H$, which we assume to be odd, and let $\tau\in H(l)$, then
\[\iota\left(W\left(k,Z_{k_1/k,\mu,\tau}\right)\right)^{\frac{l-1}{2}\frac{n}{o(\tau)}}\subseteq \rt(k_1,k,\mu).\]
\end{prop}

\begin{proof}
Let $l$ be a prime dividing $n$, let $\tau\in H(l)$ and let $x\in W(k,Z_{k_1/k,\mu,\tau})$. It follows from Lemma \ref{multiplicationmetacyclic} and Lemma \ref{metaconstructivelemma} that $\iota(x)^{\beta_l}$ is in $\rt(k_1,k,\mu_{\tilde H})$, where:
\[\begin{split}\beta_l&=\gcd\left((l-1)\frac{n}{l},(o(\tau)-1)\frac{n}{o(\tau)},\frac{3(l-1)}{2}\frac{n}{l}\right)
\\&=\gcd\left((o(\tau)-1)\frac{n}{o(\tau)},\frac{l-1}{2}\frac{n}{l}\right).
\end{split}\]
Clearly $\beta_l$ divides $\frac{l-1}{2}\frac{n}{o(\tau)}$ and so we conclude that
\[\iota(x)^{\frac{l-1}{2}\frac{n}{o(\tau)}}\in\rt(k_1,k,\mu_{\tilde H}).\]
\end{proof}

The next proposition is the main result we want to prove in this section.
\begin{prop}\label{abelcyclicresultoneinclusion}
Let $k$ be a number field and let $\G$ be a finite group such that for any class $x\in\rt(k,\G)$ there exists a tame $\G$-extension $k_1$ with Steinitz class $x$ and such that every subextension of $k_1/k$ is ramified at some primes which are unramified in $k(\zeta_{n_1})/k$.

Let $H=C(n_1)\times\dots\times C(n_r)$ be an abelian group of odd order $n$ and let $\mu$ be an action of $\G$ on $H$. We assume that the exact sequence
\[0\to H\stackrel\varphi\longrightarrow G\stackrel\psi\longrightarrow \G\to 0,\]
in which the induced action of $\G$ on $H$ is $\mu$, determines the group $G$, up to isomorphism. Then
\[\rt(k,G)\supseteq\rt(k,\G)^n \prod_{l|n}\prod_{\tau\in H(l)}W\left(k,E_{k,\mu,\tau}\right)^{\frac{l-1}{2}\frac{mn}{o(\tau)}}\]
where $E_{k,\mu,\tau}$ is the fixed field of $G_{k,\mu,\tau}$ in $k(\zeta_{o(\tau)})$, 
\[G_{k,\mu,\tau}=\left\{g\in\gal(k(\zeta_{o(\tau)})/k):\ \exists g_1\in\G,\ \mu(g_1)(\tau)=\tau^{\nu_{k,\tau}(g)}\right\}\]
and $g(\zeta_{o(\tau)})=\zeta_{o(\tau)}^{\nu_{k,\tau}(g)}$ for any $g\in \gal(k(\zeta_{o(\tau)})/k)$.
\end{prop}

\begin{proof}
Let $x\in\rt(k,\G)$ and let $k_1$ be a tame $\G$-extension of $k$, with Steinitz class $x$, and such that every subextension of $k_1/k$ is ramified at some primes which are unramified in $k(\zeta_{n_1})/k$. Thus it follows also that $k_1\cap k(\zeta_{n_1})=k$. 

By Proposition \ref{stintermediateextension}, Lemma \ref{inclusioniWEWZ}, Lemma \ref{multiplicationmetacyclic} and Proposition \ref{metaconstructive} we obtain
\[\rt(k,G)\supseteq x^n \prod_{l|n}\prod_{\tau\in H(l)}W\left(k,E_{k,\mu,\tau}\right)^{\frac{l-1}{2}\frac{mn}{o(\tau)}}.\]
We can conclude since the above inclusion holds for any $x\in\rt(k,\G)$.
\end{proof}

In this section we have only proved one inclusion concerning $\rt(k,G)$. To prove the opposite one we will need some more restrictive hypotheses. However the following lemma is true in the most general setting.

\begin{lemma}\label{primoramificatosopra}
Let $(K,k_1,k)$ be a tame $\mu$-extension and let $\P$ be a prime in $k_1$ ramifying in $K/k_1$ and let $\p$ be the corresponding prime in $k$. Then
\[x\in W(k,Z_{k_1/k,\mu,\tau})\subseteq W(k,E_{k,\mu,\tau})\subseteq \bigcap_{l|e_\P}W(k,E_{k,\mu,\tau(l)})\]
where $x$ is the class of $\p$ and $\tau$ generates $([U_\P],K/k_1)$.
\end{lemma}

\begin{proof}
Let $e_\P$ be the ramification index of $\P$ in $K/k_1$ and let $f_\p$ be the inertia degree of $\p$ in $k_1/k$. By Lemma \ref{Wram}, $\P\in H_{k_1(\zeta_{e_\P})/k_1}^{e_\P\cdot p_\infty}$ and, since the extension is tame, $\P\nmid e_\P$, i.e. $\P$ is unramified in $k_1(\zeta_{e_\P})/k_1$. Hence, by Theorem IV.8.4 in \cite{Neukirch}, $\P$ splits completely in $k_1(\zeta_{e_\P})/k_1$. It follows that the inertia degree of $\p$ in $k_1(\zeta_{e_\P})/k$ is exactly the same as in $k_1/k$, i.e. $f_\p$.

Let $u_\P\in U_\P$ be such that its class modulo $\P$ is a generator $g_\P$ of $\kappa_\P^*=U_\P/U_\P^1$. By Theorem III.8.10 and Proposition IV.6.6 in \cite{Neukirch}, $\tau=(g_\P,K/k_1)$ is an element of order $e_\P$ in $H$. 
An element $\delta\in\gal(k_1(\zeta_{e_\P})/k)$ in the decomposition group of a prime $\tilde\P$ in $k_1(\zeta_{e_\P})$ dividing $\P$, induces an automorphism of $\kappa_\P^*=\kappa_{\tilde\P}^*$ (the equality holds since $\P$ splits completely in $k_1(\zeta_{e_\P})/k_1$), given by
\[\delta(g_{\P})=g_{\P}^{\lambda_{\P,\delta}},\]
where $\lambda_{\P,\delta}$ is an integer. Thus $\zeta_{e_\P}^{\tilde\nu_{k,\mu,\tau}(\delta)}=\delta\left(\zeta_{e_\P}\right)\equiv \zeta_{e_\P}^{\lambda_{\P,\delta}} \pmod{\tilde\P}$ and, recalling that the powers of $\zeta_{e_\P}$ are distinct modulo $\tilde\P$ (since $\tilde\P\nmid e_\P$), we deduce that $\lambda_{\P,\delta}\equiv \tilde\nu_{k,\mu,\tau}(\delta)\pmod{e_\P}$. 
Recalling Proposition II.3.3 in \cite{Neukirch},
\[\begin{split}\tilde\mu_{k,\mu,\tau}(\delta)(\tau)&=\mu(\delta|_{k_1})(\tau)=\left(\delta(g_\P),K/k_1\right)\\&=\left(g_\P^{\lambda_{\P,\delta}},K/k_1\right)=\tau^{\lambda_{\P,\delta}}=\tau^{\tilde\nu_{k,\mu,\tau}(\delta)}.\end{split}\]
Thus $\delta\in \tilde G_{k_1/k,\mu,\tau}=\gal(k_1(\zeta_{e_\P})/Z_{k_1/k,\mu,\tau})$. Hence we conclude that $\p$ has inertia degree $1$ in $Z_{k_1/k,\mu,\tau}/k$ and thus it is the norm of a prime ideal in $Z_{k_1/k,\mu,\tau}$, i.e., by Proposition \ref{I.2.2Endo}, its class is in $W(k,Z_{k_1/k,\mu,\tau})$.

The proof of the inclusions
\[W(k,Z_{k_1/k,\mu,\tau})\subseteq W(k,E_{k,\mu,\tau})\subseteq W(k,E_{k,\mu,\tau(l)})\]
is trivial, using Lemma \ref{inclusioniWEWZ} and the fact that $E_{k,\mu,\tau}\supseteq E_{k,\mu,\tau(l)}$.
\end{proof}

\end{subsection}

\begin{subsection}{Realizable classes for $A'$-groups of odd order}

The next definition is technical; it will be used to make an induction argument over the order of $G$ possible. \begin{defn}
We will call a finite group $G$ \emph{good} if the following properties are verified:
\begin{enumerate}
\item For any number field $k$, $\rt(k,G)$ is a group.
\item For any tame $G$-extension $K/k$ of number fields there exists an element $\alpha_{K/k}\in k$ such that:
\begin{enumerate}
\item[(a)] If $G$ is of even order with a cyclic $2$-Sylow subgroup, then a square root of $\alpha_{K/k}$ generates the quadratic subextension of $K/k$; if $G$ either has odd order or has a noncyclic $2$-Sylow subgroup, then $\alpha_{K/k}=1$.
\item[(b)] For any prime $\p$, with ramification index $e_\p$ in $K/k$, the ideal class of
\[\left(\p^{(e_\p-1)\frac{m}{e_\p}-v_\p(\alpha)}\right)^\frac{1}{2}\]
is in $\rt(k,G)$. 
\end{enumerate}
\item For any tame $G$-extension $K/k$ of number fields, for any prime ideal $\p$ of $k$ and any rational prime $l$ dividing its ramification index $e_\p$, the class of the ideal
\[\p^{(l-1)\frac{m}{e_\p(l)}}\]
is in $\rt(k,G)$ and, if $2$ divides $(l-1)\frac{m}{e_\p(l)}$, the class of
\[\p^{\frac{l-1}{2}\frac{m}{e_\p(l)}}\]
is in $\rt(k,G)$.
\item $G$ is such that for any number field $k$, for any class $x\in\rt(k,G)$ and any integer $n$, there exists a tame $G$-extension $K$ with Steinitz class $x$ and such that every non trivial subextension of $K/k$ is ramified at some primes which are unramified in $k(\zeta_{n})/k$.
\end{enumerate}
\end{defn}

Our aim is to prove that $A'$-groups of odd order are good; but first of all at this point we need the following easy lemma.

\begin{lemma}\label{mcdle}
For any $e|m$ the greatest common divisor, for $l|e$, of the integers $(l-1)\frac{m}{e(l)}$ divides $(e-1)\frac{m}{e}$.
\end{lemma}

\begin{proof}
Let $I$ be the $\Z$-ideal generated by the integers $l-1$, for all the primes $l|e$. Then $e\equiv 1\pmod I$, since it is the product of prime factors, each one congruent to $1$ modulo $I$. Hence $e-1$ is a multiple of the greatest common divisor of the integers $l-1$ for $l|e$.

In particular for any prime $l\nmid e$, there exists an $l_1|e$, such that the power of $l$ dividing $(l_1-1)\frac{m}{e(l_1)}$ divides also $(e-1)\frac{m}{e}$.

Finally for any $l|e$ the power of $l$ dividing $(l-1)\frac{m}{e(l)}$ divides $(e-1)\frac{m}{e}$.
\end{proof}

\begin{lemma}\label{metaciclicoanalitico2}
Let $\G$ be a good group, let $H$ be an abelian group of odd order prime to that of $\G$ and let $\mu$ be an action of $\G$ on $H$.
Suppose $(K,k_1,k)$ is tamely ramified and of type $\mu$.
Let $e_\p$ be the ramification index of a prime $\p$ in $k_1/k$ and $e_\P$ be the ramification index of a prime $\P$ of $k_1$ dividing $\p$ in $K/k_1$. Then the class of
\[\left(\p^{(e_\p e_\P-1)\frac{mn}{e_\p e_\P}-v_\p(\alpha_{k_1/k}^n)}\right)^\frac{1}{2},\]
is in
\[ \rt(k,\G)^n\cdot\prod_{l|n}\prod_{\tau\in H(l)\setminus\{1\}}W\left(k,E_{k,\mu,\tau}\right)^{\frac{l-1}{2}\frac{mn}{o(\tau)}}.\]
\end{lemma}

\begin{proof}
Clearly
\[(e_\p e_\P-1)\frac{mn}{e_\p e_\P}=(e_\p-1)\frac{mn}{e_\p}+(e_\P-1)\frac{mn}{e_\p e_\P}\]
is divisible by
\[\gcd\left((e_\p-1)\frac{mn}{e_\p},(e_\P-1)\frac{mn}{e_\p e_\P}\right)\]
and, since $(m,n)=1$, i.e. also $(e_\p,e_\P)=1$, this coincides with
\[\gcd\left((e_\p-1)\frac{mn}{e_\p},(e_\P-1)\frac{mn}{e_\P}\right).\]
Thus, recalling Lemma \ref{mcdle}, 
\[\p^{(e_\p e_\P-1)\frac{mn}{e_\p e_\P}}= \p^{a_\p(e_\p-1)\frac{mn}{e_\p}+a_\P(e_\P-1)\frac{mn}{e_\P}} = \p^{a_\p(e_\p-1)\frac{mn}{e_\p}}\prod_{l|e_\P} \p^{b_{\p,l}(l-1)\frac{mn}{e_\P(l)}}. 
\]

If $\G$ either has odd order or has a noncyclic $2$-Sylow subgroup, i.e. $\alpha_{k_1/k}=1$, then we conclude by the hypothesis that $\G$ is good, by Lemma \ref{primoramificatosopra} and by the fact that any prime dividing $e_\P$ is odd.

We now assume that $\G$ is of even order with a cyclic $2$-Sylow subgroup. Again using Lemma \ref{mcdle} we can find some $c_{\p,l}$ such that
\[\begin{split}&\p^{(e_\p e_\P-1)\frac{mn}{e_\p e_\P}-v_\p(\alpha_{k_1/k}^n)}=\p^{a_\p(e_\p-1)\frac{mn}{e_\p}-v_\p(\alpha_{k_1/k}^n)}\prod_{l|e_\P} \p^{b_{\p,l}(l-1)\frac{mn}{e_\P(l)}}\\
&\qquad\qquad=\left(\p^{(e_\p-1)\frac{m}{e_\p}(a_\p-1)}\right)^n\left(\p^{(e_\p-1)\frac{m}{e_\p}-v_\p(\alpha_{k_1/k})}\right)^n\prod_{l|e_\P} \p^{b_{\p,l}(l-1)\frac{mn}{e_\P(l)}}\\
&\qquad\qquad=\prod_{l|e_\p}\p^{c_{\p,l}(l-1)\frac{mn}{e_\p(l)}(a_\p -1)}\left(\p^{(e_\p-1)\frac{m}{e_\p}-v_\p(\alpha_{k_1/k})}\right)^n\prod_{l|e_\P} \p^{b_{\p,l}(l-1)\frac{mn}{e_\P(l)}}.
\end{split}\]
We know that $\p^{(e_\p e_\P-1)\frac{mn}{e_\p e_\P}-v_\p(\alpha_{k_1/k}^n)}$ and $\p^{(e_\p-1)\frac{m}{e_\p}-v_\p(\alpha_{k_1/k})}$ are squares of ideals and that any $l$ dividing $e_\P$ is odd. It follows that $c_{\p,2}\frac{mn}{e_\p(2)}(a_\p-1)$ is even, since all the other exponents are. Recalling the hypothesis that $\G$ is good, we conclude that the class of
\[\left(\p^{(e_\p e_\P-1)\frac{mn}{e_\p e_\P}-v_\p(\alpha_{k_1/k}^n)}\right)^\frac{1}{2}\]
is in
\[\rt(k,\G)^n\cdot\prod_{l|n}\prod_{\tau\in H(l)\setminus\{1\}}W\left(k,E_{k,\mu,\tau}\right)^{\frac{l-1}{2}\frac{mn}{o(\tau)}}.\]
\end{proof}

\begin{lemma}\label{metaciclicoanalitico1}
Under the same hypotheses as in the preceding lemma, if $l|e_\p e_\P$, the class of
\[\p^{(l-1)\frac{mn}{e_\p(l)}}\]
is in
\[\rt(k,\G)^{n}\prod_{\tau\in H(l)\setminus\{1\}}W\left(k,E_{k,\mu,\tau}\right)^{\frac{l-1}{2}\frac{mn}{o(\tau)}}.\]
and, if $2$ divides $(l-1)\frac{mn}{e_\p(l)}$, the class of
\[\p^{\frac{l-1}{2}\frac{mn}{e_\p(l)}}\]
is in
\[\rt(k,\G)^{n}\prod_{\tau\in H(l)\setminus\{1\}}W\left(k,E_{k,\mu,\tau}\right)^{\frac{l-1}{2}\frac{mn}{o(\tau)}}.\]
\end{lemma}

\begin{proof}
If $l$ divides $e_\p$ the result is an obvious consequence of the fact that $\G$ is good.
For $l|e_\P$ we conclude by Lemma \ref{primoramificatosopra}. 
\end{proof}

Now we can prove the following theorem.

\begin{teo}\label{abelcyclicresult}
Let $k$ be a number field and let $\G$ be a good group.

Let $H=C(n_1)\times\dots\times C(n_r)$ be an abelian group of odd order prime to that of $\G$ and let $\mu$ be an action of $\G$ on $H$. Then
\[\rt(k,H\rtimes_\mu \G)=\rt(k,\G)^n \prod_{l|n}\prod_{\tau\in H(l)\setminus\{1\}}W\left(k,E_{k,\mu,\tau}\right)^{\frac{l-1}{2}\frac{mn}{o(\tau)}}\]
where $E_{k,\mu,\tau}$ is the fixed field of $G_{k,\mu,\tau}$ in $k(\zeta_{o(\tau)})$, 
\[G_{k,\mu,\tau}=\left\{g\in\gal(k(\zeta_{o(\tau)})/k):\ \exists g_1\in\G,\ \mu(g_1)(\tau)=\tau^{\nu_{k,\tau}(g)}\right\}\]
and $g(\zeta_{o(\tau)})=\zeta_{o(\tau)}^{\nu_{k,\tau}(g)}$ for any $g\in \gal(k(\zeta_{o(\tau)})/k)$. Furthermore $G=H\rtimes_\mu\G$ is good.
\end{teo}

\begin{proof}
Let $x\in \rt(k,H\rtimes_\mu \G)$; then $x$ is the Steinitz class of a tame extension $(K,k_1,k)$ of type $\mu$ and it is the class of a product of elements of the form
\[\left(\p^{(e_\p e_\P-1)\frac{mn}{e_\p e_\P}-v_\p(\alpha_{k_1/k}^n)}\right)^\frac{1}{2}.\]
Hence it is contained in
\[ \rt(k,\G)^n\cdot\prod_{l|n}\prod_{\tau\in H(l)\setminus\{1\}}W\left(k,E_{k,\mu,\tau}\right)^{\frac{l-1}{2}\frac{mn}{o(\tau)}}\]
by Lemma \ref{metaciclicoanalitico2} and the fact that the last expression is a group. Hence
\[\rt(k,H\rtimes_\mu \G)\subseteq  \rt(k,\G)^n\cdot\prod_{l|n}\prod_{\tau\in H(l)\setminus\{1\}}W\left(k,E_{k,\mu,\tau}\right)^{\frac{l-1}{2}\frac{mn}{o(\tau)}}.\]
The opposite inclusion is given by Theorem \ref{stdisc} and Proposition \ref{abelcyclicresultoneinclusion}.

We now show that $H\rtimes_\mu \G$ is a good group.
\begin{enumerate}
\item The first point of the definition of good groups is clear by what we have just proved about $\rt(k,H\rtimes_\mu \G)$.
\item This follows from Lemma \ref{metaciclicoanalitico2}, choosing $\alpha_{K/k}=\alpha_{k_1/k}^n$ for any extension $(K,k_1,k)$ of type $\mu$.
\item This follows from Lemma \ref{metaciclicoanalitico1}.
\item This comes from Proposition \ref{abelcyclicresultoneinclusion}.
\end{enumerate}
\end{proof}

Now we will consider direct products of good groups. We again need two lemmas.

\begin{lemma}\label{prodottodirettobuono2}
Let $\G_1$ and $\G_2$ be good groups of orders $m$ and $n$ respectively. Let us assume that $m$ and $n$ are not both even or that $\G_1$ and $\G_2$ have both non-cyclic $2$-Sylow subgroups. Let $K/k$ be a tame $\G_1\times\G_2$-extension of number fields, where $K=k_1k_2$ and $k_i/k$ are $\G_i$-extensions, let $e_\p$ be the ramification index of a prime $\p$ in $K/k$, and let
\[\alpha_{K/k}=\begin{cases}\alpha_{k_1/k}^n&\text{if $\G_1$ has even order and cyclic $2$-Sylow subgroups}\\ \alpha_{k_2/k}^m&\text{if $\G_2$ has even order and cyclic $2$-Sylow subgroups}\\1&\text{else}.\end{cases}\]
Then the class of the ideal
\[\left(\p^{(e_\p-1)\frac{mn}{e_\p}-v_\p(\alpha_{K/k})}\right)^\frac{1}{2},\]
is in
\[\rt(k,\G_1)^n\rt(k,\G_2)^m.\]
\end{lemma}

\begin{proof}
Let $\p$ be a prime ramifying in $K/k$. Let $(g_1,g_2)$ be a generator of its inertia group (it is cyclic since the ramification is tame); then $g_1$ generates the inertia group of $\p$ in $k_1/k$ and $g_2$ in $k_2/k$. Let $e_{\p,i}$ be the ramification index of $\p$ in $k_i/k$; then $e_\p=\mathrm{lcm}(e_{\p,1},e_{\p,2})$. In particular for any prime $l$ dividing $e_\p$, $e_\p(l)=\max\{e_{\p,1}(l),e_{\p,2}(l)\}$.

Let us first consider the case in which the order of $\G_1\times\G_2$ is odd or its $2$-Sylow subgroups are not cyclic. In this case $\alpha_{K/k}=1$ and, recalling Lemma \ref{mcdle},  we have
\[\begin{split}\p^{(e_\p-1)\frac{mn}{e_\p}}&=\prod_{l|e_\p}\p^{a_l(l-1)\frac{mn}{e_\p(l)}}
\\&=\prod_{\substack{l|e_\p\\e_p(l)=e_{\p,1}(l)}}\left(\p^{a_l(l-1)\frac{m}{e_{\p,1}(l)}}\right)^{n}\prod_{\substack{l|e_\p\\e_p(l)\neq e_{\p,1}(l)}}\left(\p^{a_l(l-1)\frac{n}{e_{\p,2}(l)}}\right)^{m}, \end{split}\]
where all the exponents $a_l(l-1)\frac{m}{e_{\p,1}(l)}$ and $a_l(l-1)\frac{m}{e_{\p,2}(l)}$ are clearly even. Thus, since $\G_1$ and $\G_2$ are good, the class of $\p^{\frac{1}{2}(e_\p-1)\frac{mn}{e_\p}}$ is in $\rt(k,\G_1)^n\rt(k,\G_2)^m$.

Let us now assume that $\G_1\times\G_2$ is of even order with cyclic $2$-Sylow subgroups. Thus we may suppose that the order of $\G_1$ is even, that $\G_1$ has cyclic $2$-Sylow subgroups and that the order of $\G_2$ is odd. Then

\[\begin{split}\p^{(e_\p-1)\frac{mn}{e_\p}-v_\p(\alpha_{K/k})}&=\p^{n\left((e_{\p,1}-1)\frac{m}{e_{\p,1}}-v_\p(\alpha_{k_1/k})\right)}\p^{(e_\p-1)\frac{mn}{e_\p}-(e_{\p,1}-1)\frac{mn}{e_{\p,1}}}\end{split}\]
and, recalling Theorem \ref{stdisc}, we deduce that \[\p^{(e_\p-1)\frac{mn}{e_\p}-(e_{\p,1}-1)\frac{mn}{e_{\p,1}}}\]
is the square of an ideal and we have
\[\begin{split}&\p^{(e_\p-1)\frac{mn}{e_\p}-(e_{\p,1}-1)\frac{mn}{e_{\p,1}}}\\
&\quad=\prod_{l|e_\p}\p^{a_l(l-1)\frac{mn}{e_\p(l)}}\prod_{l|e_{\p,1}}\p^{-b_l(l-1)\frac{mn}{e_{\p,1}(l)}}\\
&\quad=\prod_{\substack{l|e_\p\\e_\p(l)=e_{\p,2}(l)}}\p^{a_l(l-1)\frac{mn}{e_{\p,2}
(l)}}\prod_{\substack{l|e_{\p}\\e_\p(l)\neq e_{\p,2}(l)}}\p^{(a_l-b_l)(l-1)\frac{mn}{e_{\p,1}(l)}}\prod_{\substack{l|e_{\p,1}\\e_\p(l)= e_{\p,2}(l)}}\p^{-b_l(l-1)\frac{mn}{e_{\p,1}(l)}}.
\end{split}\]

For odd primes $l$ all the exponents in the above expression are even; we deduce that this must be true also for the component corresponding to $l=2$ (if $2|e_\p$), i.e. for $(a_2-b_2)\frac{mn}{e_{\p,1}(2)}$, and hence also for $(a_2-b_2)\frac{m}{e_{\p,1}(2)}$ since $n$ is odd.

Thus by the hypothesis that $\G_1$ and $\G_2$ are good, we easily obtain that the class of the ideal
\[\left(\p^{(e_\p-1)\frac{mn}{e_\p}-(e_{\p,1}-1)\frac{mn}{e_{\p,1}}}\right)^\frac{1}{2}\]
is in $\rt(k,\G_1)^n\rt(k,\G_2)^m$.

Now we can conclude that the class of
\[\begin{split}\left(\p^{(e_\p-1)\frac{mn}{e_\p}-v_\p(\alpha_{K/k})}\right)^\frac{1}{2}&=\p^{\frac{n}{2}\left((e_{\p,1}-1)\frac{m}{e_{\p,1}}-v_\p(\alpha_{k_1/k})\right)}\left(\p^{(e_\p-1)\frac{mn}{e_\p}-(e_{\p,1}-1)\frac{mn}{e_{\p,1}}}\right)^\frac{1}{2}\end{split}\]
is in $\rt(k,\G_1)^n\rt(k,\G_2)^m$, since also
\[\p^{\frac{n}{2}\left((e_{\p,1}-1)\frac{m}{e_{\p,1}}-v_\p(\alpha_{k_1/k})\right)}\]
is in $\rt(k,\G_1)^n$ and $\rt(k,\G_1)$ and $\rt(k,\G_2)$ are groups.
\end{proof}

\begin{lemma}\label{prodottodirettobuono1}
Under the same hypotheses as in the preceding lemma, if $l|e_\p$, the class of
\[\p^{(l-1)\frac{mn}{e_\p(l)}}\]
is in $\rt(k,\G_1)^{n}\rt(k,\G_2)^{m}$ and, if $2$ divides $(l-1)\frac{mn}{e_\p(l)}$, the class of the ideal,
\[\p^{\frac{l-1}{2}\frac{mn}{e_\p(l)}}\]
is in $\rt(k,\G_1)^n\rt(k,\G_2)^m$.
\end{lemma}

\begin{proof}
Let $l|e_\p$ and let us assume that $e_\p(l)=e_{\p,1}(l)$. 
Then
\[\p^{(l-1)\frac{mn}{e_\p(l)}}=\left(\p^{(l-1)\frac{m}{e_{\p,1}(l)}}\right)^n\]
and its class is in $\rt(k,\G_1)^{n}$, by the hypothesis that $\G_1$ is good. 
If $(l-1)\frac{mn}{e_\p(l)}$ is even then $2$ divides $(l-1)\frac{m}{e_\p(l)}$ (if $l=2$ then this is true because $2|e_\p(2)=e_{\p,1}(2)| m$ and thus by hypothesis $n$ is odd or the $2$-Sylow subgroup of $\G_1$ is not cyclic, i.e. $2$ divides $m/e_{\p,1}(2)=m/e_\p(2)$). Then
\[\p^{\frac{l-1}{2}\frac{mn}{e_\p(l)}}=\left(\p^{\frac{l-1}{2}\frac{m}{e_{\p,1}(l)}}\right)^n\]
is in $\rt(k,\G_1)^{n}$ by the assumption that $\G_1$ is good. The case $e_\p(l)=e_{\p,2}(l)$ is identical.
\end{proof}

\begin{teo}\label{directproduct}
Let $\G_1$ and $\G_2$ be good groups of orders $m$ and $n$ respectively and let us assume that $m$ and $n$ are not both even or that $\G_1$ and $\G_2$ have both non-cyclic $2$-Sylow subgroups. Then
\[\rt(k,\G_1\times\G_2)=\rt(k,\G_1)^n \rt(k,\G_2)^m.\]
Furthermore the group $\G_1\times\G_2$ is good.
\end{teo}

\begin{proof}
One inclusion is quite straightforward considering the composition of $\G_1$- and $\G_2$-extensions of $k$ with appropriate Steinitz classes and using Proposition \ref{stintermediateextension}.

The opposite inclusion follows by Lemma \ref{prodottodirettobuono2} and Theorem \ref{stdisc}.

Now again by Lemma \ref{prodottodirettobuono2} and by Lemma \ref{prodottodirettobuono1} it follows that $\G_1\times \G_2$ is good.
\end{proof}

So we obtain our most important result.

\begin{teo}
Every $A'$-group $G$ of odd order is good. In particular for any such group and any number field $k$, $\rt(k,G)$ is a subgroup of the ideal class group of $k$.
\end{teo}
\begin{proof}
Inductively, by Theorem \ref{abelcyclicresult} and Theorem \ref{directproduct}, since the trivial group is obviously good.
\end{proof}

Of course the above arguments can be used to calculate $\rt(k,G)$ explicitly for a given number field and a given $A'$-group of odd order.

Now we recall the following well-known lemma.
\begin{lemma}\label{tame2}
Let $k$ be a number field and let $\alpha\in \oo_k$ be such that $\alpha\equiv 1\pmod{4\oo_k}$. Then the extension $k(\sqrt\alpha)/k$ is tame.
\end{lemma}

\begin{proof}
By an easy calculation, $\frac{\sqrt{\alpha}+1}{2}$ is an integer, so it is in $\oo_{k(\sqrt{\alpha})}$. Now
\[\disc_{k(\sqrt\alpha)/k}\left(\left\langle 1,\frac{\sqrt{\alpha}+1}{2}\right\rangle\right)=(\alpha)\]
and so
\[\disc(k(\sqrt\alpha)/k)|(\alpha).\]
In particular it follows that $2\nmid \disc(k(\sqrt\alpha)/k)$, i.e. $2$ does not ramify in $k(\sqrt\alpha)/k$ and so the extension is tame.
\end{proof}

\begin{prop} \label{Z2banale}
Let $k$ be any number field, then
\[\rt(k,C(2))=\mathrm{cl}(k).\]
Further $C(2)$ is a good group.
\end{prop}

\begin{proof}
Let $x\in\mathrm{cl}(k)$ be any ideal class and let $\q_1$ and $\q_2$ be prime ideals in it, which are in the same ray class modulo $4$. Thanks to Proposition \ref{Tschebotarev2}, we can choose a prime ideal $\q_0$ in the ray class modulo $4$, which is inverse to that of $\q_1$ and $\q_2$.

By construction, $\q_0^2\q_1\q_2$ is principal generated by an $\alpha\equiv1\pmod{4}$.
It follows from Theorem \ref{stdisc} that 
\[D=\frac{\disc(k(\sqrt{\alpha})/k)}{\alpha}\]
is the square of a fractional ideal and by Lemma \ref{tame2} the extension $k(\sqrt{\alpha})/k$ is tame. In particular all the primes dividing $\disc(k(\sqrt{\alpha})/k)$ appear with exponent $1$ in its factorization. Then, since $(\alpha)=\q_0^2\q_1\q_2$, the only possibility for $D$ to be a square, is that it equals $\q_0^{-2}$. Then, again by Theorem \ref{stdisc}, the Steinitz class of $k(\sqrt{\alpha})/k$ is $x$.

Now it is trivial to see that $C(2)$ is good.
\end{proof}

Finally we can also use Theorem \ref{abelcyclicresult} to prove the following result about dihedral groups.

\begin{teo}
Every dihedral group $D_n$ with odd $n$ is good. In particular for any number field $k$, $\rt(k,D_n)$ is a subgroup of the ideal class group of $k$.
\end{teo}

\begin{proof}
Immediate by Theorem \ref{abelcyclicresult} and Proposition \ref{Z2banale}.
\end{proof}

\end{subsection}

\end{section}

\nocite{McCulloh}
\nocite{Long0}
\nocite{Long2}
\nocite{Endo}
\nocite{Carter}
\nocite{Massy}
\nocite{Sodaigui1}
\nocite{Sodaigui2}
\nocite{Soverchia}
\nocite{GodinSodaigui_A4}
\nocite{GodinSodaigui_ottaedri}
\nocite{ByottGreitherSodaigui} 
\nocite{CarterSodaigui_quaternionigeneralizzati}
\nocite{BrucheSodaigui}

\bibliography{bibsteinitz}
\addcontentsline{toc}{section}{Bibliography}
\bibliographystyle{abbrv}

\end{document}